\documentclass[a4paper,fleqn]{cas-sc}

\usepackage[numbers]{natbib}
\graphicspath{{./figures/}}
\usepackage{amsmath,amssymb,bm}
\usepackage{xcolor}
\usepackage{caption}
\usepackage{subcaption}
\usepackage{float}

\def\tsc#1{\csdef{#1}{\textsc{\lowercase{#1}}\xspace}}
\tsc{WGM}
\tsc{QE}

\begin{document}
\let\WriteBookmarks\relax
\def\floatpagepagefraction{1}
\def\textpagefraction{.001}

\shorttitle{}    

\shortauthors{}  

\title [mode = title]{Maximum Stress Minimization via Data-Driven Multifidelity Topology Design}  



%

\author[1]{Misato Kato}
\author[1]{Taisei Kii}
\author[1]{Kentaro Yaji}[orcid=0000-0002-4309-2043]
\cormark[1]
\ead{yaji@mech.eng.osaka-u.ac.jp}
\author[1]{Kikuo Fujita}







\affiliation[1]{organization={Department of Mechanical Engineering, Osaka University},
            addressline={2-1, Yamadaoka}, 
            city={Suita},
            postcode={565-0871}, 
            state={Osaka},
            country={Japan}}

\cortext[cor1]{Corresponding author}



\begin{abstract}
  The maximum stress minimization problem is among the most important topics for structural design. The conventional gradient-based topology optimization methods require transforming the original problem into a pseudo-problem by relaxation techniques.
  Since their parameters significantly influence optimization, accurately solving the maximum stress minimization problem without using relaxation techniques is expected to achieve extreme performance.
  This paper focuses on this challenge and investigates whether designs with more avoided stress concentrations can be obtained by solving the original maximum stress minimization problem without relaxation techniques, compared to the solutions obtained by gradient-based topology optimization.
  We employ data-driven multifidelity topology design (MFTD), a gradient-free topology optimization based on evolutionary algorithms. 
  The basic framework involves generating candidate solutions by solving a low-fidelity optimization problem, evaluating these solutions through high-fidelity forward analysis, and iteratively updating them using a deep generative model without sensitivity analysis.
  In this study, data-driven MFTD incorporates the optimized designs obtained by solving a gradient-based topology optimization problem with the $p$-norm stress measure in the initial solutions and solves the original maximum stress minimization problem based on a high-fidelity analysis with a body-fitted mesh.
  We demonstrate the effectiveness of our proposed approach through the benchmark of L-bracket.
  As a result of solving the original maximum stress minimization problem with data-driven MFTD, 
  a volume reduction of up to 22.6\% was achieved under the same maximum stress value, compared to the initial solution.
\end{abstract}



\begin{keywords}
  Topology optimization \sep Data-driven design \sep Maximum stress minimization
\end{keywords}

\maketitle

\section{Introduction}
\label{sec1}
Topology optimization, which aims to derive the best structure by optimizing material distribution within the design domain to maximize performance, is increasingly adopted in various industrial applications~\cite{Bendsoe1988,Bendsoe2004,SIGMUND2020}. While mean compliance minimization is predominant, stress-based topology optimization is actively researched for engineering applications~\cite{duysinx1998, YANG2018, KUNDU2022101716}. 

Several methods have been proposed to effectively solve the maximum stress minimization problem. However, conventional methods still face numerical and practical challenges.
As for the numerical challenge, several relaxation techniques are needed to effectively solve the maximum stress minimization and constraint problems using gradient-based methods.
This is because stress-based topology optimization typically faces three challenges: singularity, strong nonlinearity, and stress localization~\cite{Le2010}.
Singularity problems arise in density-based topology optimization, which prevent nonlinear programming algorithms from searching for the optimal solution and cause convergence to the local optima~\cite{Kirsch1990, CHENG1992, Rozvany2001, Verbart2016, Norato2022}.
To avoid this phenomenon, several techniques have been proposed to relax the stress values, such as $\varepsilon$-relaxation methods~\cite{Cheng1997} and $qp$-relaxation methods~\cite{Bruggi2008, MOON2013}.
Next, the maximum stress minimization problem is highly nonlinear and has many locally optimal solutions due to the multimodality of the solution space, making it difficult to search globally using the gradient-based methods.
Furthermore, there are challenges related to the local nature of stress evaluation points.
The computational burden increases with stress evaluation points in each element, which must be reduced by stress aggregation functions such as the $p$-norm and the Kreisselmeier-Steinhauser (KS) functions~\cite{Duysinx1998-book, Yang1996}.
Therefore, to solve the maximum stress minimization problem using gradient-based methods, it is necessary to transform the original problem into a pseudo-problem using various relaxation techniques.
These relaxation techniques cannot accurately capture the stress behavior, and the optimization results are highly sensitive to their parameters.
Next, from a practical perspective, the intermediate state of the design variable, i.e., grayscale is a significant concern. 
Gradient-based methods require that all design variables must be relaxed with continuous variables through relaxation techniques. Since grayscale is an intermediate material between solid and void, it is difficult to interpret in engineering terms and the boundaries of the optimized results are ambiguous. It must generally be removed in the design process.
In addition, density-based topology optimization typically uses a structured mesh to compute the objective function. 
It results in staircase-like boundaries in the optimized structure~\cite{LIU2016, Svard2015}, which can lead to stress concentrations at these edges.
As a result, after optimization, the designer must perform post-processing steps such as binarization and smoothing, and re-analysis using a body-fitted mesh. These steps can significantly degrade performance and cause issues such as stress concentrations, often requiring considerable time investment to achieve manufacturable designs with desirable performance.
It is necessary to achieve optimization based on stress analysis using 0/1 design variables and a body-fitted mesh.

To radically overcome the above challenges, gradient-free optimization is a promising option to deal with the 0/1 design variable field. 
Topology optimization methods using evolutionary algorithms (EAs) such as genetic algorithms~\cite{goldberg} are the representative gradient-free method~\cite{coello2007} and have been proposed in several research~\cite{Guirguis2020, WANG2005, Madeira2010}. They can perform global search even in strongly nonlinear problems. 
However, EA-based topology optimization generally increases the computational burden with an increase in design variables, namely, the \textit{curse of dimensionality}, and requires a large number of function calls of the forward analysis.
In general, only a few hundred design variables can be handled, which is very small compared to typical cases in gradient-based topology optimization.

As a new approach to achieve gradient-free optimization without compromising the design freedom of topology optimization, Yaji et al.~\cite{YAJI2022} proposed a gradient-free topology optimization framework called \textit{data-driven multifidelity topology design (MFTD)}.
The framework combines multifidelity design guided by topology optimization~\cite{Yaji2020} with data-driven topology design~\cite{Yamasaki2021}, where the solutions are updated based on EAs.
The fundamental concept is that candidate designs are generated by solving low-fidelity topology optimization problems, their objective functions are evaluated by high-fidelity forward analysis, and iteratively updated by a deep generative model that corresponds to crossover in EAs until the desired solutions are obtained.
Data-driven MFTD has been shown to be applicable to topology optimization problems that are difficult to solve directly with conventional gradient-based methods, such as minimax~\cite{kato2023, kii2024} and turbulent flow problems~\cite{YAJI2022}.
The maximum stress minimization problem has been addressed by Yamasaki et al.~\cite{Yamasaki2021}. However, there have not been sufficient discussions of the specific challenges inherent in this problem and the effectiveness of the framework in addressing those challenges. 

Therefore, this paper focuses on the challenges of the maximum stress minimization problem faced by conventional methods and investigates whether solutions with more avoided stress concentrations can be obtained by solving the original maximum stress minimization problem based on data-driven MFTD, compared to the solutions obtained by conventional methods.
In this paper, "fidelity" is defined as the accuracy of the maximum stress value.
A low-fidelity problem is a pseudo-problem that incorporates relaxation techniques such as intermediate density, approximate processing, and a structured mesh, which are typically handled by conventional gradient-based methods. On the other hand, a high-fidelity problem is an original problem that deals with the maximum value of the true stress obtained by forward analysis using black-and-white design and a body-fitted mesh.
Specifically, we first generate several design solutions by gradient-based topology optimization using the $p$-norm, as a low-fidelity optimization. These solutions are then used as the initial solutions of the proposed framework, and they are iteratively updated based on a high-fidelity evaluation using an original maximum stress minimization problem to improve their performance.
To improve the convergence of data-driven MFTD, this paper incorporates latent crossover, a method recently proposed by Kii et al. to efficiently perform crossover in the latent space of ~\cite{kii2024}.
Through numerical examples, we discuss the challenges of the conventional method and the effectiveness of data-driven MFTD for the maximum stress minimization problem.

\begin{figure*}[t]
  \centering\includegraphics[width=\linewidth]{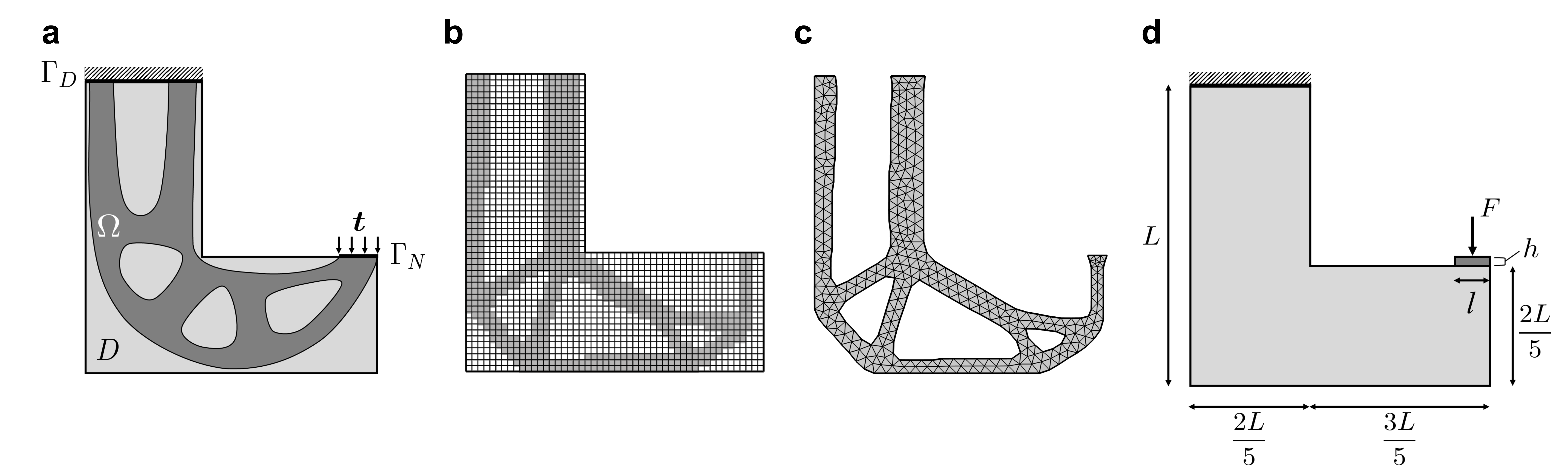}
  \caption{Problem settings for L-bracket. (a) Design domain $D$ and the shape $\Omega$; (b) Discretization with structured mesh; (c) Discretization with body-fitted mesh; (d) Boundary conditions and dimensions for numerical examples in Section \ref{sec5}.}\label{design_domain}
\end{figure*}
\section{Problem settings}
\label{sec2}
In gradient-based topology optimization, the maximum stress minimization problem is usually defined as a volume-constrained single-objective optimization problem.
Since EA-based optimization methods can naturally deal with multi-objective optimization problems, we consider the bi-objective optimization of maximum stress minimization and volume minimization as the original problem.

Fig.~\ref{design_domain}a shows the L-bracket widely used as a benchmark for maximum stress minimization problems.
$\Omega$ and $D\backslash\Omega$ denote the material and void domains in the design domain $D$, respectively.
The material domain $\Omega$ is designed inside a pre-fixed design space $D$. The structure is fixed on $\Gamma_D$ and $\bm{t}$ is the surface force applied on $\Gamma_N$.
Since the initial structure concentrates stress at the re-entrant corner, the material in this area is aggressively removed during optimization to reduce the maximum stress.

Let us consider the continuous form of a maximum stress minimization problem.
A formulation of the bi-objective optimization problem can be defined by
\begin{equation}\label{opt1}
  \begin{aligned}
  &\text{find}&& \Omega \subseteq D,\\
  &\text{that minimize}&&\sigma_\text{max} = \max\{\sigma_\text{vm}(\bm{x}) \mid \forall \bm{x} \in \Omega\},\\
  & &&V = \int_{\Omega}d\Omega.\\
  \end{aligned}
\end{equation}
Herein, $\sigma_\text{vm}$ is the von Mises stress and $\bm{x}$ is the point inside the material domain $\Omega$.
Note that Eq.~(\ref{opt1}) deals with the maximum value of the true stress.
Since the maximum stress is not differentiable, density-based methods cannot solve this original problem directly and must transform it into a pseudo-problem using relaxation techniques.
The conventional formulation using the relaxation techniques is described in Section~\ref{sec3}.

The analysis domain is typically discretized using the finite element method to calculate the evaluation functions.
In density-based optimization, a structured mesh is commonly used as shown in Fig.~\ref{design_domain}b. However, since the optimized structure depends on the shape of the mesh and stress concentration occurs at the edges, it is not suitable for the evaluation of maximum stress values. In this paper, we calculate the evaluation function in Eq.~(\ref{opt1}) using a body-fitted mesh as shown in Fig.~\ref{design_domain}c, which is necessary to more accurately evaluate the maximum stress.

\section{Conventional approach}\label{sec3}
In gradient-based topology optimization, it is necessary to transform the original problem into a pseudo-problem using various relaxation techniques to solve the maximum stress minimization problem. In this section, we introduce general relaxation techniques and formulation using these techniques.

We discuss the general density-based topology optimization for solving the maximum stress minimization problem. 
Based on a finite element analysis using a structured mesh and the solid isotropic material with penalization (SIMP) \cite{CHENG1992, Rozvany2001, Cheng1997, Rozvany1994}, we defined a pseudo-density $\rho\in[0,1]$ at each element as the design variable.

\subsection{Density Filter}\label{sec:df}
To avoid a checkerboard pattern, the SIMP method introduces the density filter~\cite{bourdin2001, BRUNS2001} for the design variable $\rho_i$ as follows:
\begin{equation}\label{eq_H}
  \begin{aligned}
  \tilde{\rho}_i &= \frac{\sum_{j \in \Omega_i}w_j\rho_i}{\sum_{j \in \Omega_i}w_j}, \quad i = 1,2,...,n,\\
  \end{aligned}
\end{equation}
where $\Omega_i$ is in the influence domain of the element $i$. The weight factor $w_j$ is calculated at each element $j$ inside $\Omega_i$ as follows:
\begin{equation}\label{eq_w}
  \begin{aligned}
  w_j &= \frac{r_0 - r_j}{r_0},\\
  \end{aligned}
\end{equation}
where $r_0$ is the filter radius, and $r_j$ is the distance between the center points of each elements $i$ and $j$.
The density filter can control the minimum length scale, and avoid mesh-dependence and stress concentration due to small and sharp structures.

\subsection{Equilibrium Equation}
Under the assumption of linear elasticity and static behavior, the discrete form of the equilibrium equation can be represented as $\bm{Ku} = \bm{F}$.
$\bm{K}$ can be adequately built by using the stiffness matrix of an element $\bm{K}_{i}$ and the modified Young's modulus $E_{\text{SIMP}}$, given by:
\begin{equation}\label{eq_K}
\begin{aligned}
\bm{K} &= \sum_{i = 1}^{n}E_{\text{SIMP}}(\tilde{\rho}_i)\bm{K}_{i}\\
E_{\text{SIMP}}(\tilde{\rho}_i) &= E_{\min} + \tilde{\rho}^{p}_i(E_0 - E_{\min}),\\
\end{aligned}
\end{equation}
where $E_0$ and $E_{\min}$ are Young's modulus of the solid $(\tilde{\rho_i} = 1)$ and void $(\tilde{\rho_i} = 0)$ phases, respectively. 
$E_{\min}$ is a small positive real numbers to avoid numerical instabilities.
And $p$ is the penalization parameter to promote binarization to achieve black-and-white designs finally. Typically, $p = 3$ is often used in topology optimization problems concerning the maximum stress minimization or constraint, as with the mean compliance~\cite{duysinx1998, Le2010, Holmberg2013}.

\subsection{Stress Relaxation}
It is widely known that maximum stress minimization problems face a singularity problem, in which nonlinear programming algorithms cannot reach degenerate regions of the design space that often contain global optimum, and converges to a local solution~\cite{Le2010}. Specifically, stress at the void could increase rapidly as one or more of the design variables tends to be zero.
Hence, stress relaxation methods that smooth the design space are necessary to avoid singularity and stabilize the optimization process.
Various stress relaxation methods have been proposed, e.g., the $\varepsilon$-relaxation and smooth envelope functions \cite{Cheng1997,Rozvany1992}. In this study, we use the $qp$-parametrization \cite{Le2010, Bruggi2008, Holmberg2013}, one of the general relaxation methods.

The vector of stress at evaluation point $i$ can be written in Voigt notation as
\begin{equation}\label{eq_stress}
\bm{\sigma}_{i} = (\sigma_{ix}\,\sigma_{iy}\,\sigma_{iz}\,\tau_{ixy}\,\tau_{iyz}\,\tau_{izx})^\mathsf{T}.\\
\end{equation}
The penalized and relaxed stress measure $\hat{\bm{\sigma}}_{i}$ interpolating stress values for intermediate density is given as
\begin{equation}\label{eq_relax}
\begin{aligned}
\hat{\bm{\sigma}}_{i}(\tilde{\rho}_i) &= \eta(\tilde{\rho}_i)\bm{\sigma}_{i},\\
\eta(\tilde{\rho}_i) &= \tilde{\rho}_i^q,\\
\end{aligned}
\end{equation}
where $q$ is the penalization parameter, typically $q = 0.5$ is used~\cite{Le2010, Holmberg2013}. This method also has the effect of penalizing intermediate values of the material density and promoting binarization. In addition, the following holds for both extreme values:
\begin{align}\label{eq_relax_01}
\hat{\bm{\sigma}}_{i}(\tilde{\rho}_i = 1) &= \bm{\sigma}_{i},\\
\hat{\bm{\sigma}}_{i}(\tilde{\rho}_i = 0) &=\lim_{\tilde{\rho}_i \to 0}\hat{\bm{\sigma}}_{i}(\tilde{\rho}_i) = 0.
\label{eq_relax_01_2}
\end{align}
Herein, Eq.~(\ref{eq_relax_01_2}) justifies that the $qp$-parametrization can avoid the singularity problem \cite{Kocvara2012}.

\subsection{$P$-norm Stress Measure}
In this study, the relaxed von Mises stress of Eq.~(\ref{eq_relax}) is used as the stress measure in the optimization procedure.
In the optimization problem~(\ref{opt1}), we replace the objective function expressed as the maximum stress with $\sigma_{\max} = \underset{i}{\max}(\hat{\sigma}_{\text{vm},i})$.

The maximum stress is not differentiable, so typically it needs to be approximated using a global function in gradient-based methods. We use $p$-norm stress measure given by
\begin{equation}\label{eq_p}
\begin{aligned}
\sigma_\text{PN} &= \left(\sum_{i=1}^{n} \hat{\sigma}_{\text{vm},i}^P \right)^{1/P}.\\
\end{aligned}
\end{equation}
Herein, the $p$-norm $\sigma_\text{PN}$ approaches the maximum stress $\sigma_{\max}$ when the parameter $P \rightarrow \infty$, whereas numerical computation becomes impossible. Therefore, a finite value is used for $P$, but the larger the value of $P$, the greater the numerical instability.
Thus, the optimization highly depends on the selection of this parameter, and it is necessary to select an appropriate value for the stress norm parameter $P$ that provides good search performance and numerical stability. In previous works, although $P = 8$ was often chosen from the viewpoint of numerical stability during optimization~\cite{Le2010}, it cannot accurately capture the maximum stress.
In this paper, we employ the continuation method~\cite{Bendsoe2004} to achieve both a large $P$ setting and numerical stability. It allows numerical stability by increasing the stress norm parameter $P$ at every pre-determined interval during optimization. The effectiveness of this method is discussed in Section~\ref{sec5}.

\subsection{Problem Formulation}
In this study, the original optimization problem of Eq.~(\ref{opt1}) has two objectives: maximum stress minimization and volume minimization. To solve this bi-objective problem using the conventional gradient-based topology optimization methods, we use the $\varepsilon$-constraint approach \cite{Marler2004, Miettinen2012} to replace the problem with a single objective optimization problem, the volume-constrained maximum stress minimization problem.
Consequently, the optimization problem can be formulated as follows:
\begin{equation}\label{opt2}
\begin{aligned}
&\text{find}&& \bm{\rho} = {\rho_i}\,(i = 1,2,...,n)\\
&\text{that minimize}&& \sigma_\text{PN} = \left(\sum_{i=1}^{n}\hat{\sigma}_{\text{vm},i}^P \right)^{1/P}\\
&\text{subject to}&&V = \sum_{i=1}^{n}v_i\rho_i \leq \overline{V} \leq V_{\max}\\
& && \rho_{i}\in[0,1],\\
\end{aligned}
\end{equation}
where $v_i$ is the element $i$ solid volume, $\overline{V}$ is the volume constraint, and $V_{\max}$ is the design domain volume.
Herein, Eq.~(\ref{opt2}), which includes several relaxation techniques, can be defined as a pseudo-problem compared to the original optimization problem of Eq.~(\ref{opt1}). 
We solve this optimization problem by using the method of moving asymptotes (MMA)~\cite{Svanberg1987}, 
one of the popular gradient-based optimizer in the research community of topology optimization.
\begin{figure*}[t]
  \centering\includegraphics[width=\linewidth]{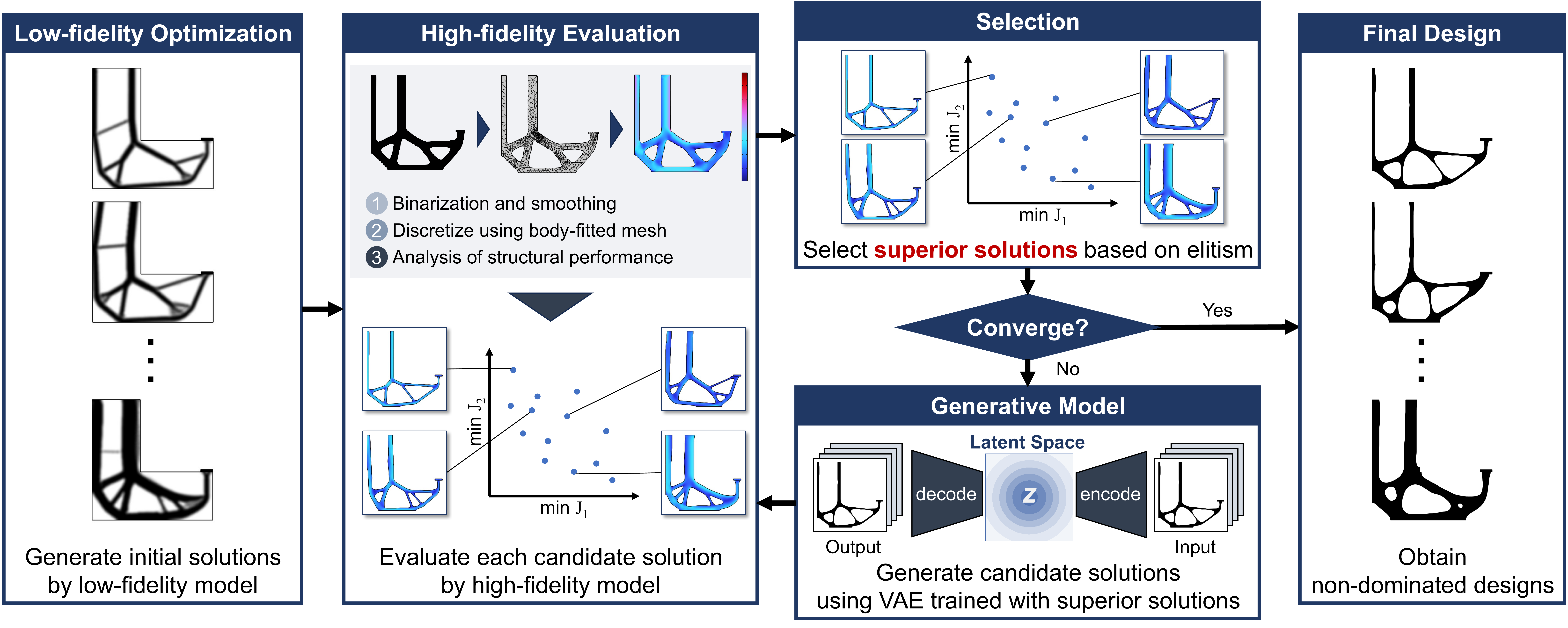}
  \caption{Schematic illustration of data-driven multifidelity topology design}\label{framework}
\end{figure*}

\section{Proposed approach}\label{sec4}
In this study, we tackle the original maximum stress minimization problem~(\ref{opt1}) based on data-driven multifidelity topology design (MFTD)~\cite{YAJI2022} that is a gradient-free topology optimization framework under a high degree of design freedom.
Specifically, we investigates whether optimized designs obtained by solving the optimization problem (\ref{opt2}) which is handled by the conventional gradient-based method, can be improved on the original topology optimization problem (\ref{opt1}).
It should be emphasized that the proposed approach in this study is a framework for optimization that updates the solutions based on the results of high-fidelity stress analysis using a body-fitted mesh without the relaxation techniques used in gradient-based methods.

The procedures of data-driven MFTD are shown in Fig.~\ref{framework}. 
Each step is briefly described below.
Note that this paper omits a mutation-like operation proposed in the original framework for simplicity and the detailed procedures can be found in the original paper~\cite{YAJI2022}.

\subsection{Low-fidelity Optimization}
In contrast to the original problem to be solved, a pseudo-problem, i.e., a low-fidelity optimization problem, is defined.
A pseudo-optimization problem incorporating design parameters, called seeding parameters, is solved by gradient-based topology optimization to generate a variety of design candidates.
In this study, we use the relaxed problem formulation in Eq.~(\ref{opt2}) as the low-fidelity optimization problem and the seeding parameter is the volume constraint $\overline{V}$. Various patterns of initial designs for the framework are generated based on the $\varepsilon$-constraint method.

\subsection{High-fidelity Evaluation}
In this step, all the candidates are evaluated by high-fidelity model on the original objective space in Eq.~(\ref{opt1}), namely, the maximum von Mises stress $\sigma_{\max}$ and the volume fraction $V/V_{\max}$.
The high-fidelity model treated here is a model as shown in Fig.~\ref{design_domain}c, in which binarization, smoothing, and body-fitted mesh are applied to low-fidelity model represented by pixel as shown in Fig.~\ref{design_domain}b. This allows the analysis on the original problem~(\ref{opt1}).
Note that the framework only requires the forward analysis of the original high-fidelity model without any gradient information on the objective and constraint functions.

\subsection{Selection}
Based on the results of the high-fidelity evaluation, superior candidates so-called elite solution are selected from the candidate solutions using an elite strategy of EAs. In this study, the non-dominated sorting genetic algorithm II (NSGA-II), one of the representative selection algorithms, is used to rank and select candidates based on the Pareto dominance relation in the objective space \cite{Deb2022}.
Based on the EA strategy, high-fidelity evaluation, selection, and generation of candidate solutions are repeated until the Pareto front converges.

\subsection{Generative Model}
The aim in this step is to generate new candidate solutions with the characteristics of the selected elite solution in selection by using a generative model. The important point here is that the material distributions can be updated without sensitivity analysis.
We use a variational autoencoder (VAE) \cite{Kingma2013}, which is one of the representative deep generative models.  
As shown in Fig.~\ref{framework}, a VAE consists of two neural networks, namely, an encoder and a decoder. It has the ability to extract the information from the high-dimensional input data, and to compressed it into a lower-dimensional manifold called latent space.
In the standard VAE, the latent variable $\bm{z}$ is defined as
\begin{equation}
	\bm{z}=\bm{\mu}+\bm{\sigma}\circ\bm{\varepsilon},
\end{equation}
where $\bm{\mu}$ is the mean, $\bm{\sigma}$ is the variance, $\circ$ is the element-wise product, and $\bm{\varepsilon}$ is a random vector following the standard normal distribution. 
With this architecture, a VAE uses the same dataset for inputs and outputs to perform unsupervised learning to construct the latent space.
In this way, it is expected to extract essential features of the training data by compressing high-dimensional input and output data into a low-dimensional latent space.

Consequently, by sampling from this latent space and generating material distributions, it is possible to generate new candidate solutions that inherit the characteristics of elite solutions.
In this study, to reduce the randomness in sampling process, we employ latent crossover, a scheme for sampling candidates intensively from meaningful regions in the latent space~\cite{kii2024}.

\section{Results and Discussion}\label{sec5}
In this section, we present numerical examples of a maximum stress minimization problem and demonstrate the challenges of conventional gradient-based topology optimization and the effectiveness of the proposed approach with data-driven MFTD.

\subsection{Design Settings}
We deal with the L-bracket shown in Fig.~\ref{design_domain}d, which gives the dimensions, loads, and boundary conditions. All the constants are dimensionless values, and $L=2$, $l = 0.2$, $h = 0.04$. To avoid stress concentration, the total force $F = 1$ applied to the top boundary of the non-design domain. The Young's modulus is set as $E_0 = 1$ for the solid and $E_{\min} = 10^{-9}$ for the void. Poisson's ratio is $\nu = 0.3$ and the filter radius is $r_0 = 0.05$. The 2D solid element considering plane stress state is employed. For the low-fidelity model, the number of finite elements in the structured mesh is 25,600, and for high-fidelity model, the body-fitted mesh is applied to each solution.
\begin{figure*}[h!]
  \centering\includegraphics[width=\linewidth]{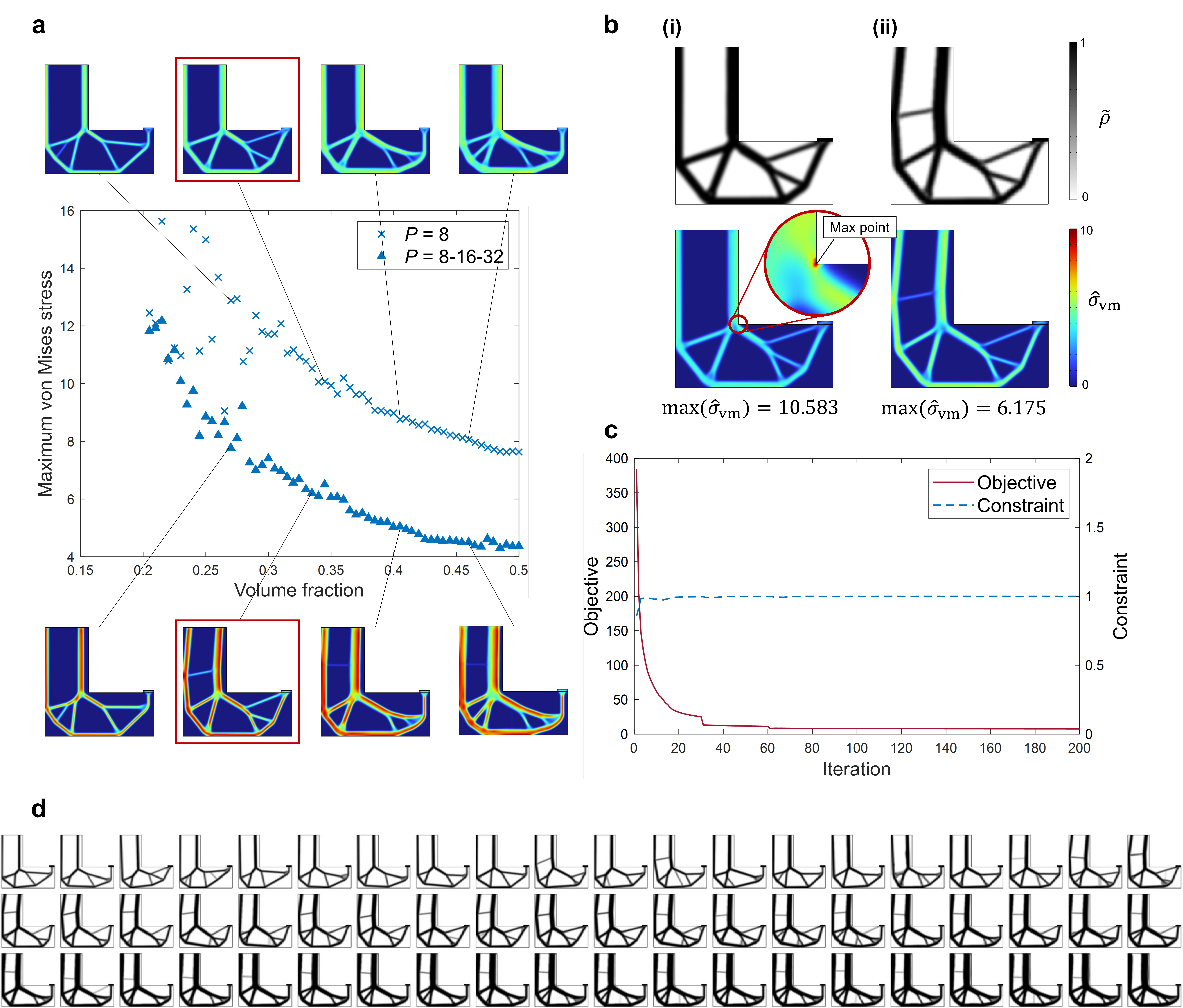}
  \caption{Effect of the continuation method on the maximum von Mises stress. (a) Objective space of the relaxed maximum von Mises stress $\max(\hat{\sigma}_{\text{vm}})$ and volume fraction $V/V_{\max}$ with the stress distributions of the solutions; (b) Optimized designs with the filtered density $\tilde{\rho}$ and their stress distributions in the same volume constraint $\overline{V}/V_{\max} = 0.335$ under (i) the fixed parameter $P = 8$ and (ii) the continuation method $P = \textrm{8-16-32}$; (c) Convergence history of the objective function and the normalized volume constraint in the low-fidelity optimization using the continuation method $P = \textrm{8-16-32}$; (d) Initial solutions for data-driven MFTD generated by the low-fidelity optimization using the continuation method $P = \textrm{8-16-32}$.}\label{low}
\end{figure*}

\subsection{Generation of Initial Solutions}\label{initial}
We generate various design candidates by adjusting the volume constraints $\overline{V}/V_{\max}$ from 0.2 to 0.5 in 0.05 increments.
The discrete adjoint method is used for sensitivity analysis to derive gradient information, and the MMA is used as the gradient-based optimizer~\cite{Svanberg1987}.
The convergence criterion is simply set to the maximum iteration of 200. The move limit of the MMA is set to 0.05. 

In order to verify whether there is room for improvement in the optimized designs obtained by the conventional method, it is necessary to prepare the best initial solutions as much as possible.
Therefore, we employ the continuation method for the parameter $P$ of $p$-norm stress measure to obtain solutions with more avoided stress concentrations.
First, we test the effectiveness of the continuation method, in which $P$ is increased every 30 steps as 8-16-32 by comparing the solutions obtained using the fixed stress norm parameter $P = 8$.

Fig.~\ref{low}a shows the effect of the continuation method on the maximum von Mises stress with the stress distributions of the solutions. Here, the maximum von Mises stress is the maximum value of the relaxed stress $\sigma_{\max} = \underset{i}{\max}(\hat{\sigma}_{\text{vm},i})$ defined as Eq.~(\ref{eq_relax}).
The results indicate that the continuation method achieves a better Pareto front and obtained the solutions with a more uniform stress distribution compared with the $P = 8$ fixed case.
Fig.~\ref{low}b shows the optimized designs and their stress distributions in the same volume constraint $\overline{V}/V_{\max} = 0.335$ under the fixed parameter $P = 8$ and the continuation method $P = \textrm{8-16-32}$, which are the same solutions shown by the red box in Fig.~\ref{low}a.
The re-entrant corner of a solution under $P = 8$ is filled with the material as shown in Fig.~\ref{low}b(i), causing stress concentrations, while the continuation method more accurately captures the maximum, hence material in this area is actively removed and stress concentrations are relaxed as shown in Fig.~\ref{low}b(ii).
Therefore, this method can avoid such a singularity problem.
Fig.~\ref{low}c shows the convergence history of the objective function and the normalized volume constraint for the solution under the continuation method $P = \textrm{8-16-32}$ shown in Fig~\ref{low}b(ii).
The objective function here is the $p$-norm stress measure defined as Eq.~(\ref{eq_p}).
As shown in Fig.~\ref{low}c, the objective and constraint functions converge well enough in all solutions, and it can be said that, at least under the investigated conditions, the gradient-based topology optimization method with the continuation method has resulted in optimized structures.
From the above, we adopt the continuation method in low-fidelity optimization and use 60 optimized designs as the initial designs for data-driven MFTD.
Fig.~\ref{low}d shows the optimized designs of them.
Here, let us focus on the structure of the obtained results. As shown in Fig.~\ref{low}d, various structures appear even when the volumes are close to each other.
At first glance, there is no clear relationship between structural changes and volume increases, and the number of members also varies. However, very similar structures can be observed even when the volumes differ.

\begin{figure*}[t!]
  \centering\includegraphics[width=\linewidth]{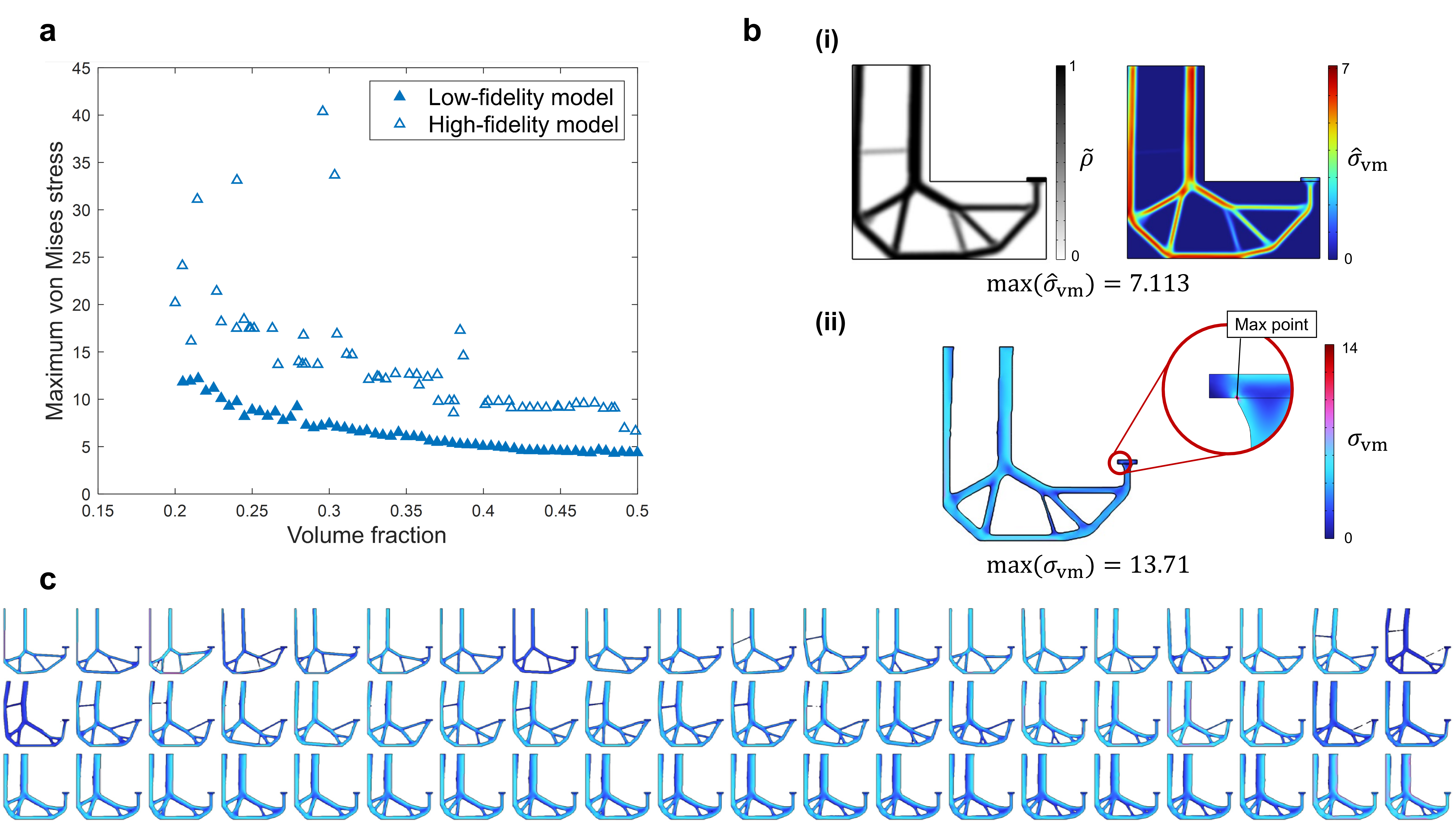}
  \caption{Effect of the fidelity on the maximum von Mises stress. (a) Objective space of the maximum von Mises stress $\sigma_{\max}$ and volume fraction $V/V_{\max}$, in which each stress value indicates the relaxed stress $\max(\hat{\sigma}_{\text{vm}})$ (low-fidelity model) and the true stress $\max(\sigma_{\text{vm}})$ (high-fidelity model), respectively; (b) Optimized design and its stress distribution in different analysis model: (i) the low-fidelity model with the filtered density $\tilde{\rho}$, (ii) the high-fidelity model; (c) Initial solutions analyzed with the high-fidelity model.}\label{fidelity}
\end{figure*}

Next, Fig.~\ref{fidelity} shows the effect of the fidelity of the analysis model on the maximum von Mises stress. 
The result for the low-fidelity model shows the maximum relaxed stress of the solutions obtained by low-fidelity optimization, which is the same as the result for the continuation method $P = \textrm{8-16-32}$ in Fig.~\ref{low}a.
On the other hand, the result for the high-fidelity model shows the maximum values of the true stress obtained through analyzing them after post-processing and discretization with a body-fitted mesh.
As shown in Fig.~\ref{fidelity}a, the maximum stress exhibits significant variability.
This is due to stress concentrations resulting from structural changes caused by post-processing, in addition to the difference in stress indices.
Fig.~\ref{fidelity}b shows the optimized design and its stress distribution in different analysis model, the low-fidelity and high-fidelity model.
Post-processing has resulted in the disappearance of members composed of intermediate densities, and the structure has also changed.
In the low-fidelity model, stress is almost uniformly distributed throughout the whole structure, whereas in the high-fidelity model, stresses are concentrated locally at a single point in the loading area, and the maximum stress values are different.
Fig.~\ref{fidelity}c shows the stress distribution of all designs. It can be noticed that stresses are concentrated locally in most of the solutions and some designs have disconnected members.
It should be emphasized that this is a serious problem in actual design because it is difficult to predict the variations in stress values due to the fidelity of the analysis model.
\begin{figure*}[t!]
  \centering\includegraphics[width=\linewidth]{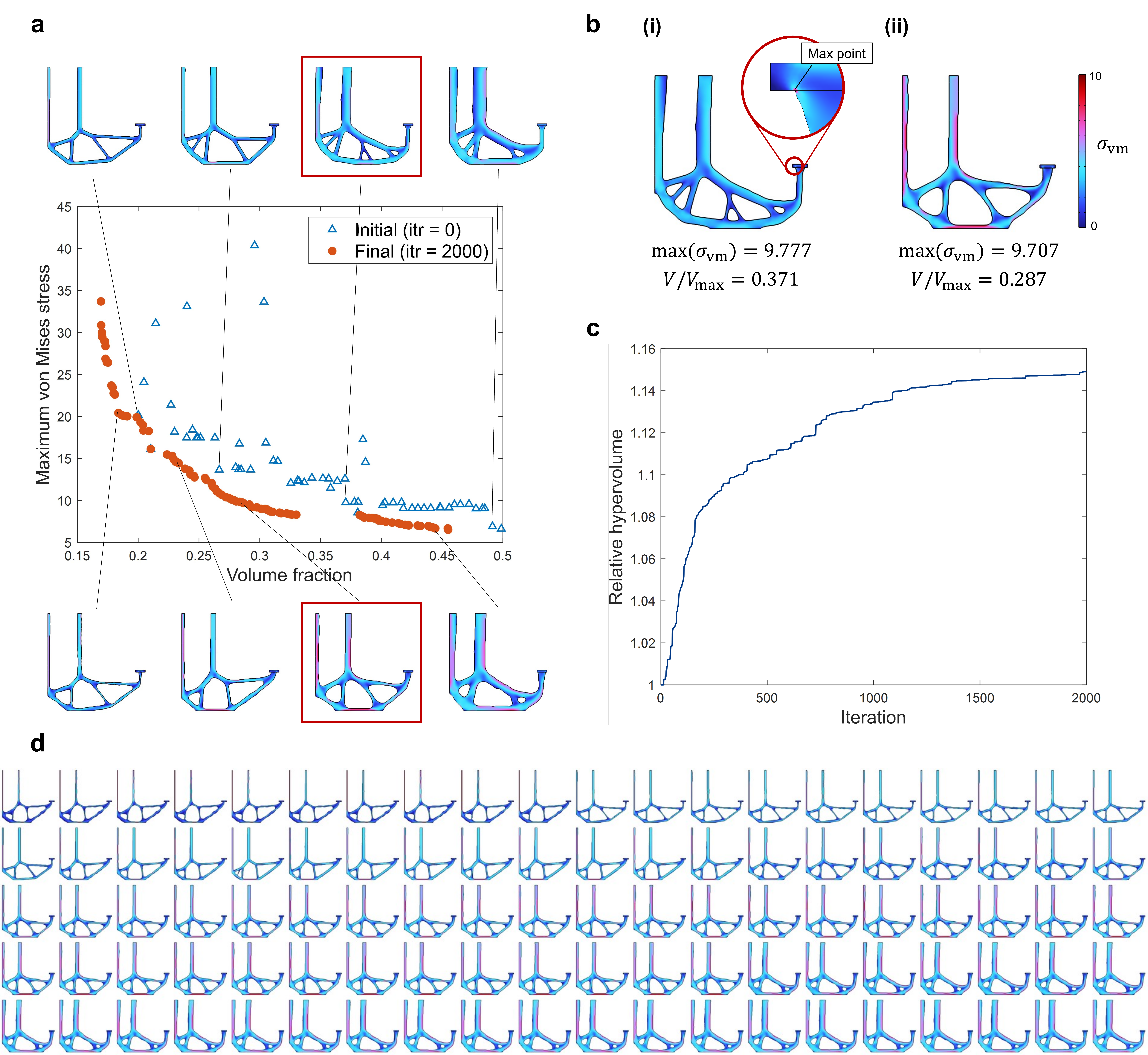}
  \caption{Optimization results with data-driven MFTD. (a) Objective space with initial and final designs; (b)(i) An initial solution and (ii) A final solution under the similar maximum von Mises stress; (c) Convergence history of the hypervolume indicator; (d) Final solutions by deta-driven MFTD.}\label{high}
\end{figure*}

\subsection{Improvement of Pareto Solutions}
This study investigates the potential enhancement of performance in designs optimized through the gradient-based topology optimization method by applying further optimization using data-driven MFTD. It aims to demonstrate the efficacy of the gradient-free approach for the maximum stress minimization problem.

The initial solutions obtained in Fig.~\ref{fidelity}c is input to data-driven MFTD to produce the Pareto solutions. Then, we tackle to solve the original optimization problem (\ref{opt1}) based on high-fidelity stress analysis. The convergence criterion is the maximum iteration number of 2,000, after which the solutions were considered to be converged. Here, At each optimization step, 100 elite solutions survive and are trained by a VAE, after which 100 new candidate solutions are generated. 
The architecture of the VAE is a simple neural network based on a multilayer perceptron with one hidden layer each for the encoder and decoder. In this study, the number of neurons in the input/output and hidden layers is 25,600 and 512, respectively. For the latent space, the mean $\bm{\mu}$, the variance $\bm{\sigma}$ and the latent variable $\bm{z}$ are structured on 8 neurons.

Fig.~\ref{high}a shows the objective space with the stress distribution of some initial and final solutions.
The Pareto front has been progressed from the initial solutions, indicating that the conventional gradient-based topology optimization has room for improvement in terms of solving the original maximum stress minimization problem. 
Fig.~\ref{high}b shows the initial and final designs with almost the same maximum stress value, which are the same solutions shown by the red box in Fig.~\ref{high}a.
The objective function values are $\max(\sigma_\text{vm}) = 9.777$ with $V/V_{\max} = 0.371$, and $\max(\sigma_\text{vm}) = 9.707$ with $V/V_{\max} = 0.287$, respectively.
Specifically, a volume reduction of 22.6\% was achieved under almost the same stress value, compared to the initial solution.
Focusing on the structures, it can be observed that the final solution has fewer members and a simpler structure than the initial solution. Additionally, the solutions obtained by conventional gradient-based topology optimization are linear and uniform in member thickness. On the other hand, since there is no limit on the thickness of the members in data-driven MFTD, the final solution features more rounded holes and re-entrant corner to relax stress compared to the initial solution.
Furthermore, focusing on stress distributions, the initial solution has a localized concentration of stress, whereas the final solution has a more uniform stress distribution due to stress dispersion.

Fig.~\ref{high}c shows the convergence history of the hypervolume indicator which is a measure often used to assess the performance of multi-objective evolutionary algorithms~\cite{Shang2021}. In the case of two objective functions, it is represented by the area formed by the reference point and the Pareto front in the objective space. In this study, the reference point is the maximum values of the objective functions in the initial solutions.
As shown in Fig.~\ref{high}c, the hypervolume gradually improves up to 2,000 iterations, indicating that the Pareto front continues to progress.
Note that the number of iterations is relatively large compared to other optimization problems investigated in the original paper~\cite{YAJI2022}.
Fig.~\ref{high}d shows the final designs. It can be observed that the initial solutions are various structures depending on the volume, while the final solutions have a relatively similar topology. This result indicates that although the Pareto solutions are not necessarily the globally optimal, it can be expected that the comprehensive solution search via data-driven MFTD has yielded a promising topology among numerous local solutions. It can also be said that data-driven MFTD yielded superior solutions that could not be reached by the conventional method, since structures with more avoided stress concentrations were obtained that are not typically found in the gradient-based method.

\section{Conclusion}\label{sec6}
In this paper, we focused on the remaining challenges of the maximum stress minimization problem solved by the standard gradient-based topology optimization method and proposed a new framework to overcome those challenges.
To accurately solve the maximum stress minimization problem, we focused on data-driven MFTD, a gradient-free topology optimization method, and proposed an optimization framework based on high-fidelity stress analysis using a body-fitted mesh.
This framework derived initial solutions by solving the gradient-based topology optimization using the $p$-norm stress measure, and they are updated by a VAE based on the manner of EAs, without using sensitivity analysis.
It was confirmed that solutions with more avoided stress concentrations were obtained compared to the initial solutions obtained by the gradient-based method, as the Pareto front has been progressed by data-driven MFTD. The comprehensive solution search using data-driven MFTD for the original maximum stress minimization problem resulted in structures with characteristics not commonly seen in conventional methods, achieving designs with reduced stress concentrations.
we also achieved to derive one promising topology among many local solutions.

As for future work, we consider that our framework can be applied to further complex and practical optimization problems considering strong geometrical nonlinearity such as buckling and large deformations under the original maximum stress minimization or constraint problem.

\section*{Acknowledgements}\label{Acknowledgements}
This work was supported by JSPS KAKENHI Grant Numbers 23K26018 and JP24KJ1639.









\bibliographystyle{elsarticle-num}
\bibliography{cas-refs}



\end{document}